\documentclass[11pt,a4paper]{article}
\usepackage[T1]{fontenc}
\usepackage[utf8]{inputenc}
\usepackage{amsmath,amssymb,amsthm}
\newtheorem{proposition}{Proposition}
\usepackage{newpxtext,newpxmath}
\usepackage{microtype}
\usepackage[margin=26mm,top=25mm,bottom=25mm]{geometry}
\usepackage[numbers,sort&compress]{natbib}
\usepackage{fancyhdr}
\usepackage{needspace}
\usepackage{booktabs}
\usepackage{caption}
\usepackage{graphicx}
\usepackage{placeins}
\usepackage{xurl}
\usepackage[hidelinks]{hyperref}
\hypersetup{pdftitle={Seven Coincidences and Sixty Symmetries},pdfsubject={Binomial coincidences and their symmetries},pdfauthor={Ernesto Lupercio and Evgeniya Shaydurova}}
\begin{document}
\thispagestyle{plain}
\begin{center}
{\Large\bfseries Seven Coincidences and Sixty Symmetries}\\[0.8em]
{\large Ernesto Lupercio$^{1}$ and Evgeniya Shaydurova$^{2}$}\\[0.55em]
{\small $^1$Departamento de Matem\'aticas, CINVESTAV, 07360 Ciudad de M\'exico, Mexico}\\
{\small $^2$Universidad Abierta y a Distancia de M\'exico, Mexico}\\[0.4em]
{\small Correspondence: \href{mailto:lupercio@math.cinvestav.mx}{lupercio@math.cinvestav.mx}}
\end{center}
\vspace{0.8em}
\begin{abstract}
Motivated by John Baez's question about seven exceptional coincidences among binomial coefficients, we construct seven exterior-power identities for representations of $A_5$. Their common symmetry comes from simplex edges and icosahedral vertices, axes, and faces, together with characters of stabilizers. We enumerate all ordered pairs of representation isomorphism classes in the prescribed dimensions. Apart from the pair of trivial actions, the permutation pair at $1540$ and the complex monomial pair at $7140$ are unique. At $3003$, Klein's twelve-point divisor gives a noncanonical intertwiner factored through polynomial division, differentiation, and residues.
\end{abstract}
\vspace{0.5em}
\noindent
John Baez asked for structural explanations of seven unusual equalities between binomial coefficients. The smallest compares pairs chosen from sixteen objects with triples chosen from ten:
\[
\binom{16}{2}=\frac{16\cdot15}{2}
=120=
\frac{10\cdot9\cdot8}{3\cdot2}=\binom{10}{3}.
\]
Dividing by $2$ and $6$ accounts for the orderings of a pair and a triple. A structural explanation would also relate the objects being counted. In the same essay, Baez reports Paul Schwahn's representation-theoretic explanation and asks for explanations of the other six equalities. That question led to this article \citep{Baez2026Coincidences}.

Taking complements gives an example of such a correspondence. A choice of three objects from ten determines the seven left over, as in
\[
\{1,4,9\}\quad\longleftrightarrow\quad\{2,3,5,6,7,8,10\}.
\]
Applying the operation again returns the original choice, so it pairs every three-element subset with exactly one seven-element subset and proves
\[
\binom{10}{3}=\binom{10}{7}.
\]
Two finite sets of the same cardinality can always be put in bijection by enumerating their elements; complementation is more informative because the correspondence is determined by the subsets themselves and respects every permutation of the ten labels. For the seven equalities, we shall seek linear isomorphisms compatible with a group action.

The group $A_5$ acts in all seven constructions. It is the group of sixty rotations preserving an icosahedron and, equivalently, the group of even permutations of five letters. We use its actions on simplex edges and on the vertices, axes, and faces of an icosahedron to construct representations whose exterior powers explain the seven equalities.

\FloatBarrier
\section{Repeated entries in Pascal's triangle}\label{sec:pascal}

The seven values appear together in a short addition to David Singmaster's 1971 paper, where he reported six numbers occurring six times in Pascal's triangle and a seventh, \(3003\), occurring eight times \citep[p. 386]{Singmaster1971}. Writing the binomial coefficients beside their common values gives
\begingroup
\allowdisplaybreaks[3]
\begin{align*}
\binom{16}{2}&=\binom{10}{3}=120,\\[3pt]
\binom{21}{2}&=\binom{10}{4}=210,\\[3pt]
\binom{56}{2}&=\binom{22}{3}=1540,\\[3pt]
\binom{78}{2}&=\binom{15}{5}=\binom{14}{6}=3003,\\[3pt]
\binom{120}{2}&=\binom{36}{3}=7140,\\[3pt]
\binom{153}{2}&=\binom{19}{5}=11628,\\[3pt]
\binom{221}{2}&=\binom{17}{8}=24310.
\end{align*}
\endgroup
The first expression in each row is the triangular number \(\binom m2=m(m-1)/2\). The row for \(3003\) contains three distinct pairs of indices, giving three positions in Pascal's triangle but just one of the seven values.

The triangle had a history long before Pascal's name became attached to it: Jia Xian's eleventh-century work, preserved in part through Yang Hui's account of 1261, used its coefficients in procedures for extracting roots, alongside the rule that forms an entry by adding its two neighbors above \citep{MacTutorJiaXian}. We will fix an integer and ask where it occurs in Pascal's triangle.

Pascal's work on the triangle is associated with 1654, and his treatise appeared posthumously in 1665. Its applications include combinations and the division of stakes in an interrupted game; he takes the four letters \(A,B,C,D\) and lists their six unordered pairs \citep[pp. 469--470]{Pascal1923Triangle}:
\[
AB,\quad AC,\quad AD,\quad BC,\quad BD,\quad CD.
\]
The pair \(AA\) is excluded because it repeats a letter, while \(AB\) and \(BA\) count as the same choice, so these conventions give \(\binom42=6\); every permutation of the four letters also permutes the six pairs. Adding a fifth letter will turn this elementary observation into the first of our finite symmetry constructions.

\begin{figure}[!htbp]
\centering
\includegraphics[width=0.98\linewidth]{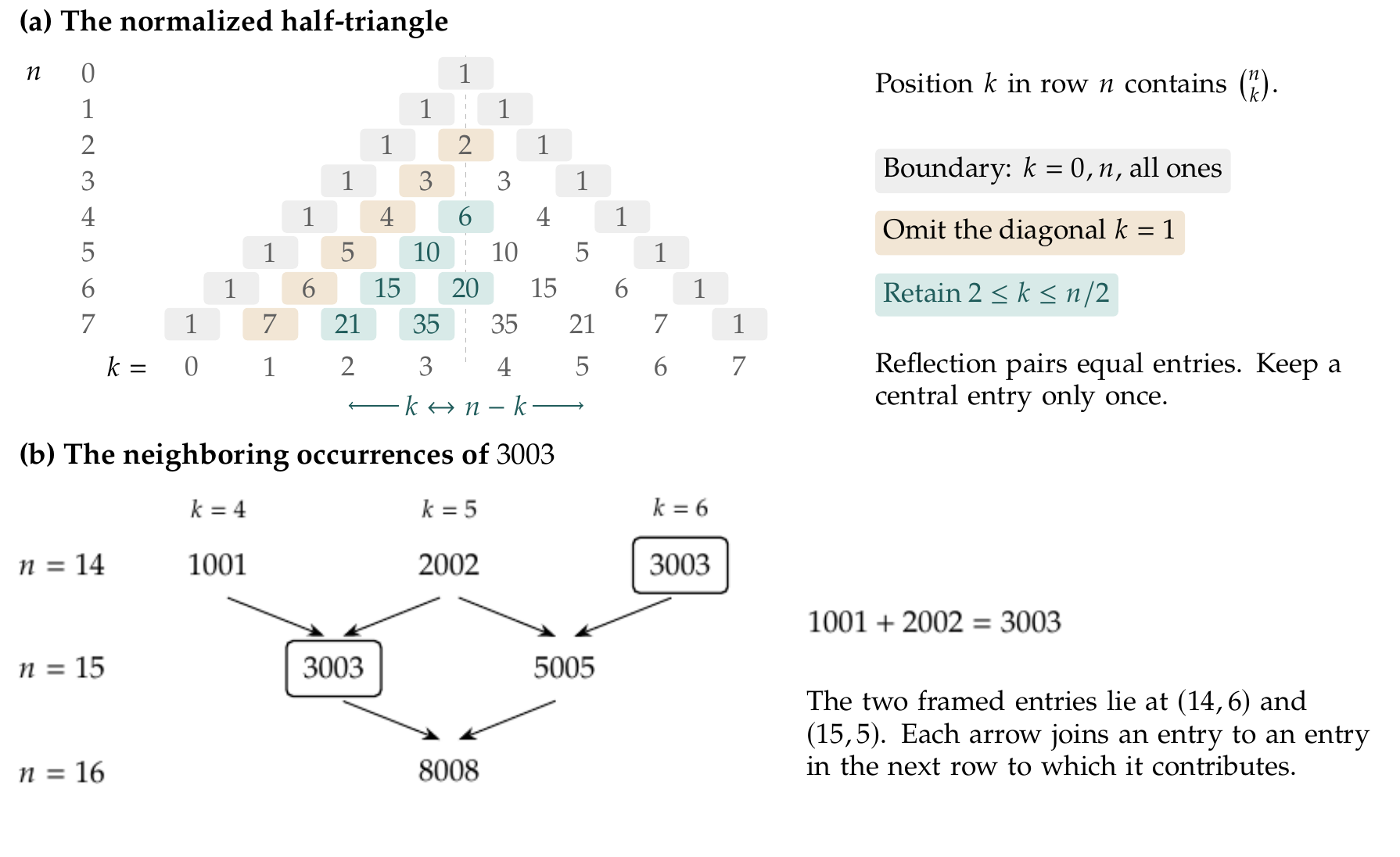}
\caption{(a) Normalization of Pascal's triangle: retain $2\leq k\leq n/2$ (shaded), excluding the boundary ones and the $k=1$ diagonal and choosing one entry from each reflected pair. (b) The two framed entries are $\binom{14}{6}$ and $\binom{15}{5}$. Arrows show the local Pascal additions; $\binom{78}{2}$ lies elsewhere.}\label{fig:pascal}
\end{figure}

Put \(\binom nk\) in position \(k\) of row \(n\), numbering both from zero. Complementation gives the reflection \(\binom nk=\binom n{n-k}\); the boundary entries are ones, and every positive integer \(a\) also occurs as \(\binom a1\). Singmaster counted all positions, writing
\[
N(a)=\#\left\{(n,k):n\geq0,\ 0\leq k\leq n,\ \binom nk=a\right\},
\]
where \(\#\) denotes cardinality. For fixed \(a>1\), this number is finite: an entry with \(1\leq k\leq n-1\) is at least \(n\), so an occurrence of \(a\) requires \(n\leq a\). Singmaster conjectured that \(N(a)\) is bounded by a constant independent of \(a>1\) \citep{Singmaster1971}. For the list of nontrivial collisions, we use the normalized half-triangle
\[
2\leq k\leq n/2.
\]
Figure~\ref{fig:pascal}(a) illustrates this convention: it removes the boundary and the neighboring diagonal \(k=1\), chooses one representative under reflection, and counts a central entry in an even row once.

In row fourteen the entries at positions four, five, and six are \(1001,2002,3003\), and Pascal's addition rule gives
\[
\binom{14}{4}+\binom{14}{5}
=1001+2002=3003
=\binom{15}{5}=\binom{14}{6}.
\]
Singmaster displayed this pattern in 1975 \citep[p. 295]{Singmaster1975}; Figure~\ref{fig:pascal}(b) shows these additions. In the normalized half-triangle, the three occurrences are \((14,6),(15,5),(78,2)\). Adding their reflections and the two occurrences \((3003,1)\) and \((3003,3002)\) gives the eight occurrences
\[
\begin{gathered}
(14,6),(14,8),\qquad (15,5),(15,10),\\
(78,2),(78,76),\qquad (3003,1),(3003,3002).
\end{gathered}
\]
Hence \(N(3003)=8\). Pascal's rule gives the equality between the entries in rows fourteen and fifteen. Their equality with $\binom{78}{2}$ is a separate identity.

Douglas Lind obtained an infinite family of neighboring-row equalities in 1968, and Singmaster treated it again in 1975 [\citealp[Section 5]{Lind1968}; \citealp[p. 296]{Singmaster1975}]. Its first nontrivial members are
\[
\binom{15}{5}=\binom{14}{6},\qquad
\binom{104}{39}=\binom{103}{40},\qquad
\binom{714}{272}=\binom{713}{273}.
\]
Their indices can be expressed through the Fibonacci numbers \(0,1,1,2,3,5,\ldots\), each new term being the sum of its two predecessors. Section~\ref{sec:arithmetic}, ``The same numbers in arithmetic, derives these indices from a quadratic equation. Here the family accounts for the neighboring equality in the \(3003\) row. The equality with $\binom{78}{2}$ does not belong to this family.

In 1997, Benjamin M. M. de Weger conjectured that this family and the seven exceptional values exhaust the nontrivial collisions \citep[Conjecture A]{deWeger1997}. This proposed classification would imply \(N(a)\leq8\), but it also specifies which equalities occur, whereas a uniform bound alone says nothing about their locations. Blokhuis, Brouwer, and de Weger verified the proposed list in two search ranges: all common values at most \(10^{60}\), and all rows through \(10^6\) \citep[Theorem 1]{BlokhuisBrouwerdeWeger2017}. These finite verifications leave the unrestricted conjecture open.

The classification problem asks which collisions occur; Baez's question asks what structure can explain these seven equalities. We construct representations of $A_5$ and compare their exterior powers.

Baez's n-Category Caf\'e post of 24 September 2026 presents combinatorial arguments for three of the equalities, developed from a suggestion of Mike Stay \citep{Baez2026Cafe}. These include
\[
\binom{78}{2}=\binom{14}{6}=3003.
\]
Our icosahedral construction gives the different pairing
\[
\binom{78}{2}=\binom{15}{5}.
\]

\FloatBarrier
\section{Choices and their symmetries}\label{sec:choices}

Let \(X\) be a finite set of labels, and assign to each label \(x\) a basis vector \(e_x\). The resulting vector space, written \(\mathbb C[X]\), consists of the formal sums \(\sum_{x\in X}a_xe_x\) with complex coefficients. A permutation \(g\) of the labels extends to a linear transformation by
\[
\rho(g)e_x=e_{g(x)},\qquad
\rho(g)\Bigl(\sum_x a_xe_x\Bigr)=\sum_x a_xe_{g(x)}.
\]
Performing two permutations successively gives \(\rho(gh)=\rho(g)\rho(h)\), with \(h\) acting first. An assignment $g\mapsto\rho(g)$ of invertible linear maps satisfying this composition rule is a \emph{representation} of the group; the maps need not permute a basis, although that is how our first examples arise. We take these permutation representations over \(\mathbb C\), including those obtained from real geometric configurations.

For three labels, the transposition \(g=(12)\), which exchanges 1 and 2 and fixes 3, has matrix
\[
\rho(g)=\begin{pmatrix}0&1&0\\1&0&0\\0&0&1\end{pmatrix},
\qquad \operatorname{Tr}\rho(g)=1.
\]
The columns record the images of the three basis vectors. A diagonal entry is one when its label is fixed, so the trace of a permutation matrix counts fixed labels. This example uses $S_3$, the permutation group on three labels.

The permutation action preserves both the constant vectors and the vectors whose coordinates sum to zero. For \(v=(x_1,\ldots,x_n)\), let \(\bar x=(x_1+\cdots+x_n)/n\). Then
\[
v=\bar x(1,\ldots,1)
   +\bigl(v-\bar x(1,\ldots,1)\bigr).
\]
The first term lies in a line fixed pointwise by every permutation, and the second has coordinate sum zero. Each subspace is \emph{invariant}, meaning that every group element maps it into itself. Their intersection is zero, and every vector has the displayed unique decomposition, so the space is their \emph{direct sum}. We write \(\mathbf1\) for a line on which every group element acts as the identity, and \(rW\) for the direct sum of \(r\) copies of \(W\). Dimensions add under direct sum.

To compare two spaces together with their actions, we require the identification to commute with every symmetry. For representations \(U,V\), a linear map \(T:U\to V\) is \emph{equivariant} when
\[
T\rho_U(g)=\rho_V(g)T\qquad\text{for every }g.
\]
Our complement correspondence provides an example: let \(U\) have basis the three-element subsets of ten labels and \(V\) the seven-element subsets. Let \(T:U\to V\) be the linear map determined by \(T(e_S)=e_{S^c}\), where \(S^c\) is the complement of \(S\) in the ten-element label set. A permutation of the ten labels satisfies
\[
T\rho_U(g)e_S=e_{(gS)^c}=e_{g(S^c)}=\rho_V(g)Te_S.
\]
Taking complements commutes with every permutation, so \(T\) is equivariant; taking complements twice gives the identity, so \(T\) is invertible. We write \(U\cong V\) when there exists an invertible equivariant linear map \(U\to V\).

Alongside the direct sum, the \emph{tensor product} \(U\otimes V\) has product basis \(e_i\otimes f_j\), one vector for each pair of basis labels; its dimension is \((\dim U)(\dim V)\), and the map \((u,v)\mapsto u\otimes v\) is bilinear. The induced action is
\[
g(u\otimes v)=(gu)\otimes(gv).
\]
When the two factors are the same, exchanging them either preserves a symmetric product or reverses the sign of an exterior product; this distinction determines whether repeated basis vectors are allowed.

In the exterior square \(\Lambda^2W\), products are written \(u\wedge v\) and satisfy
\[
u\wedge v=-v\wedge u,\qquad u\wedge u=0.
\]
For \(W=\mathbb C^3\), a basis is \(e_1\wedge e_2,e_1\wedge e_3,e_2\wedge e_3\). More generally, \(\Lambda^kW\) has basis \(e_{i_1}\wedge\cdots\wedge e_{i_k}\) with \(i_1<\cdots<i_k\), giving \(\dim\Lambda^kW=\binom{\dim W}{k}\). The wedge basis is indexed by the \(k\)-element subsets; permuting the factors multiplies a wedge by the sign of that permutation, and the action applies \(g\) to every factor before restoring increasing index order.

The transposition $g=(12)$ acts by
\[
\begin{aligned}
g(e_1\wedge e_2)&=-e_1\wedge e_2,\\
g(e_1\wedge e_3)&= e_2\wedge e_3,\\
g(e_2\wedge e_3)&= e_1\wedge e_3.
\end{aligned}
\qquad
g|_{\Lambda^2W}=\begin{pmatrix}-1&0&0\\0&0&1\\0&1&0\end{pmatrix}.
\]
The trace is $-1$ on $\Lambda^2W$ and $1$ on $W$. The representations therefore cannot be equivariantly isomorphic: an invertible equivariant map would conjugate the two action matrices, and conjugation preserves trace.

In the symmetric square \(\operatorname{Sym}^2W\), exchanging the factors leaves the product unchanged, so a repeated factor need not vanish. For a two-dimensional space the comparison is
\[
\operatorname{Sym}^2\mathbb C^2:
     \quad e_1^2,\ e_1e_2,\ e_2^2,
\qquad
\Lambda^2\mathbb C^2:\quad e_1\wedge e_2.
\]
A basis of \(\operatorname{Sym}^2W\) consists of the products \(e_i e_j\) with \(1\leq i\leq j\leq n\), giving \(n+\binom n2=n(n+1)/2\) vectors.

Write \(\chi_W(g)=\operatorname{Tr}(g|_W)\). If the eigenvalues of \(g\) on \(W\) are \(\lambda_1,\ldots,\lambda_n\), its eigenvalues on the exterior square are \(\lambda_i\lambda_j\) for \(i<j\), and
\[
\left(\sum_i\lambda_i\right)^2
=\sum_i\lambda_i^2+2\sum_{i<j}\lambda_i\lambda_j.
\]
The first sum on the right is the trace of \(g^2\). Consequently
\begin{equation}\label{eq:exterior-square}
\chi_{\Lambda^2W}(g)
=\frac{\chi_W(g)^2-\chi_W(g^2)}2.
\end{equation}
Including the diagonal products gives \(\chi_{\operatorname{Sym}^2W}(g)=(\chi_W(g)^2+\chi_W(g^2))/2\). For every exterior degree at once, choosing either 1 or \(t\lambda_i\) from each factor gives
\begin{equation}\label{eq:exterior-generating}
\sum_{r=0}^{n}\chi_{\Lambda^rW}(g)t^r
   =\prod_{i=1}^{n}(1+t\lambda_i)=\det(1+tg|_W).
\end{equation}
The coefficient of \(t^r\) is the trace on \(\Lambda^rW\). Group elements of finite order are diagonalizable over \(\mathbb C\), so these eigenvalue calculations apply to every element used below.

The trace function \(\chi_W\) is called the \emph{character}. Two finite-dimensional complex representations of a finite group are equivariantly isomorphic if and only if their characters agree at every element. A representation is \emph{irreducible} if it has no proper nonzero invariant subspace. Every finite-dimensional complex representation of a finite group is a direct sum of irreducibles, whose multiplicities are determined by character inner products. Section A.1 of the companion \citep{Companion2026} proves these facts by averaging an inner product and an operator over the group. They turn finitely many trace calculations into an existence proof for an equivariant isomorphism.

The spinorial example credited to Schwahn in Baez's account can now be read as an identity of representations \citep{Baez2026Coincidences}. The group \(\operatorname{Spin}(10)\) is a double cover of the rotations of ten-dimensional Euclidean space: it maps onto that rotation group with two elements over each rotation. Its ten-dimensional vector representation $V_{10}$ and sixteen-dimensional half-spin representation $S_{10}^{+}$ satisfy the classical identity
\[
\Lambda^2S_{10}^{+}\cong\Lambda^3V_{10},\qquad
\dim\Lambda^2S_{10}^{+}=\binom{16}{2},\quad
\dim\Lambda^3V_{10}=\binom{10}{3}.
\]
The isomorphism is a classical representation identity, verified by spin characters in Section D.1 of the companion; the equality of dimensions follows from it. We shall later realize its restriction to \(A_5\) using tensors on a four-dimensional space. For the next identity, both representations come directly from permutations.

\FloatBarrier
\section[Five letters and the identity at 210]{Five letters and the identity at $210$}\label{sec:five-letters}

Add a fifth letter to Pascal's four, and label the letters \(1,\ldots,5\). The determinant of a permutation matrix is \(+1\) or \(-1\); a permutation is called \emph{even} in the first case. Composition multiplies determinants, so the even permutations form a subgroup, \(A_5\). Multiplication by the transposition \((12)\) pairs every even permutation with an odd one, giving \(|A_5|=5!/2=60\). The ten unordered pairs are the edges of the complete graph on five vertices, or of a regular four-simplex. Write \(R_5\) for the five-letter permutation space and \(A\) for its four-dimensional zero-sum subspace, so \(R_5=\mathbf1\oplus A\). The ten edges give a permutation space \(E\). We will compare \(\Lambda^2U\) with \(\Lambda^4E\), where \(U=\mathbf1\oplus2E\) has dimension twenty-one.

Consider \(g=(12)(34)\). It fixes the edges \(12\) and \(34\), while exchanging
\[
13\leftrightarrow24,\qquad14\leftrightarrow23,\qquad
15\leftrightarrow25,\qquad35\leftrightarrow45.
\]
Here \(ij\) denotes the unordered edge \(\{i,j\}\). The two fixed edges contribute two eigenvalues \(+1\), and each exchange contributes \(+1,-1\). Hence
\[
\det(1+tg|_E)=(1+t)^2(1-t^2)^4,
\qquad\chi_E(g)=2.
\]
More generally, a cycle of length \(\ell\) has the \(\ell\)th roots of unity as eigenvalues. From \(\prod_{\lambda^\ell=1}(x-\lambda)=x^\ell-1\) we obtain the cycle factor \(1-(-t)^\ell\), so fixed labels, exchanged pairs, and longer cycles can all be handled by the same product.

Elements \(g,g'\in G\) are \emph{conjugate} if \(g'=hgh^{-1}\) for some \(h\in G\). Their action matrices are then conjugate, so their traces agree. In \(A_5\) the classes have representatives of orders \(1,2,3,5,5\), denoted by \(1,2A,3A,5A,5B\). There are fifteen double transpositions, twenty three-cycles, and twenty-four five-cycles. A five-cycle commutes with just its five powers among the permutations of five labels, so its conjugacy class in \(A_5\) has \(60/5=12\) elements; the twenty-four five-cycles therefore form two classes.

We can now check the identity at \(210\) on these five classes. For \(g=(12)(34)\), we have already found \(\chi_E(g)=2\). Since \(g\) is an involution, meaning that \(g\ne1\) and \(g^2=1\), \(\chi_U(g)=1+2\cdot2=5\), \(\chi_U(g^2)=21\). Formula~\eqref{eq:exterior-square} gives \(\chi_{\Lambda^2U}(g)=(25-21)/2=2\). On the other side, writing \([t^4]\) for the coefficient of \(t^4\), we find
\[
[t^4](1+t)^2(1-t^2)^4
=[t^4](1+2t+t^2)(1-4t^2+6t^4+\cdots)
=6-4=2.
\]
The fourth exterior power of \(E\) therefore has the same trace at \(g\).

On edges a three-cycle has one fixed edge and three cycles of length three, and a five-cycle has two cycles of length five. Hence
\[
\chi_E=(10,2,1,0,0),\qquad
\chi_U=(21,5,3,1,1),\qquad
\chi_{\Lambda^2U}=(210,2,3,0,0).
\]
For the three-cycle, \(g^2\) has the same trace as \(g\), giving \((3^2-3)/2=3\) on the left; for either five-cycle class it gives \((1^2-1)/2=0\). The right-hand generating polynomials are \((1+t)(1+t^3)^3\) and \((1+t^5)^2\), with coefficients of \(t^4\) three and zero. The dimension is \(\binom{10}{4}=210\), so the two characters agree on all five conjugacy classes and the character theorem proves
\[
\Lambda^2(\mathbf1\oplus2E)\cong\Lambda^4E,
\qquad\binom{21}{2}=\binom{10}{4}=210.
\]
The restriction to even permutations matters: a single transposition on five labels gives traces \(30\) and \(-14\) on these two spaces, obtained respectively from \((9^2-21)/2\) and \([t^4](1+t)^4(1-t^2)^3\).

For higher exterior powers, the generating polynomial gives a recurrence in traces of powers. Let \(e_r=\chi_{\Lambda^rV}(g)\) and \(p_i=\chi_V(g^i)\), so that \(p_i\) is the trace of the \(i\)th power, not the \(i\)th power of the trace. For \(E_g(t)=\prod_j(1+t\lambda_j)=\sum_r e_rt^r\),
\[
\frac{E_g'(t)}{E_g(t)}
=\sum_j\frac{\lambda_j}{1+t\lambda_j}
=p_1-p_2t+p_3t^2-\cdots.
\]
Multiplying by \(E_g(t)\) and comparing coefficients gives
\[
e_0=1,\qquad e_1=p_1,\qquad
2e_2=p_1^2-p_2,\qquad
3e_3=p_1e_2-p_2e_1+p_3;
\qquad
re_r=\sum_{i=1}^{r}(-1)^{i-1}p_i e_{r-i}.
\]
This is Newton's recurrence. Together with the square formula it computes the exterior-power characters on the five conjugacy classes of \(A_5\); the complete calculation is given in companion Section A \citep{Companion2026}.

\section[Vertices, axes, and 3003]{Vertices, axes, and $3003$}\label{sec:vertices}

The next identity uses the standard geometric identification of \(A_5\) with the rotations of an icosahedron. Its twelve vertices and six pairs of opposite vertices give permutation spaces \(Q\) and \(P\). Each set is a single \emph{orbit}: any one of its labels can be carried to any other by a rotation. The latter pairs determine six unoriented axes; an axis may be fixed even when the rotation exchanges its endpoints. Let \(B\subset P\) be the five-dimensional zero-sum subspace, an invariant complement to the constant line. Functions on the vertices split into even and odd parts under the antipodal map:
\[
f_+(v)=\frac{f(v)+f(-v)}2,\qquad
f_-(v)=\frac{f(v)-f(-v)}2.
\]
Even functions on the vertices descend to functions on the six axes, while odd functions are determined by one value on each opposite pair and form a six-dimensional space \(C\). Rotations commute with the antipodal map, so these projections are equivariant and
\[
P=\mathbf1\oplus B,\qquad Q=\mathbf1\oplus B\oplus C.
\]
We use the action \((gf)(v)=f(g^{-1}v)\), which sends a basis function supported at \(v\) to the one supported at \(gv\). Rotations of orders two, three, and five fix respectively zero, zero, and two vertices; on axes their cycle types are \(1^2 2^2,3^2,1\,5\). Counting the fixed labels and subtracting the characters of the invariant summands gives Table~\ref{tab:characters}. Set \(S=2B\oplus C\), of dimension sixteen, for use at \(11628\) and in the spinorial identity. Here and below a \emph{module} means a vector space together with its group action.

\begin{center}
\begin{minipage}{\linewidth}
\captionof{table}{Characters of the geometric modules. The identity column gives their dimensions.}\label{tab:characters}
\centering
\renewcommand{\arraystretch}{1.12}
\begin{tabular}{lrrrrr}
\toprule
Module & \(1\) & \(2A\) & \(3A\) & \(5A\) & \(5B\)\\
\midrule
\(A\) & 4 & 0 & 1 & -1 & -1\\
\(B\) & 5 & 1 & -1 & 0 & 0\\
\(C\) & 6 & -2 & 0 & 1 & 1\\
\(E\) & 10 & 2 & 1 & 0 & 0\\
\(P\) & 6 & 2 & 0 & 1 & 1\\
\(Q\) & 12 & 0 & 0 & 2 & 2\\
\(S\) & 16 & 0 & -2 & 1 & 1\\
\bottomrule
\end{tabular}
\end{minipage}
\end{center}

\begin{figure}[!htbp]
\centering
\includegraphics[width=0.98\linewidth]{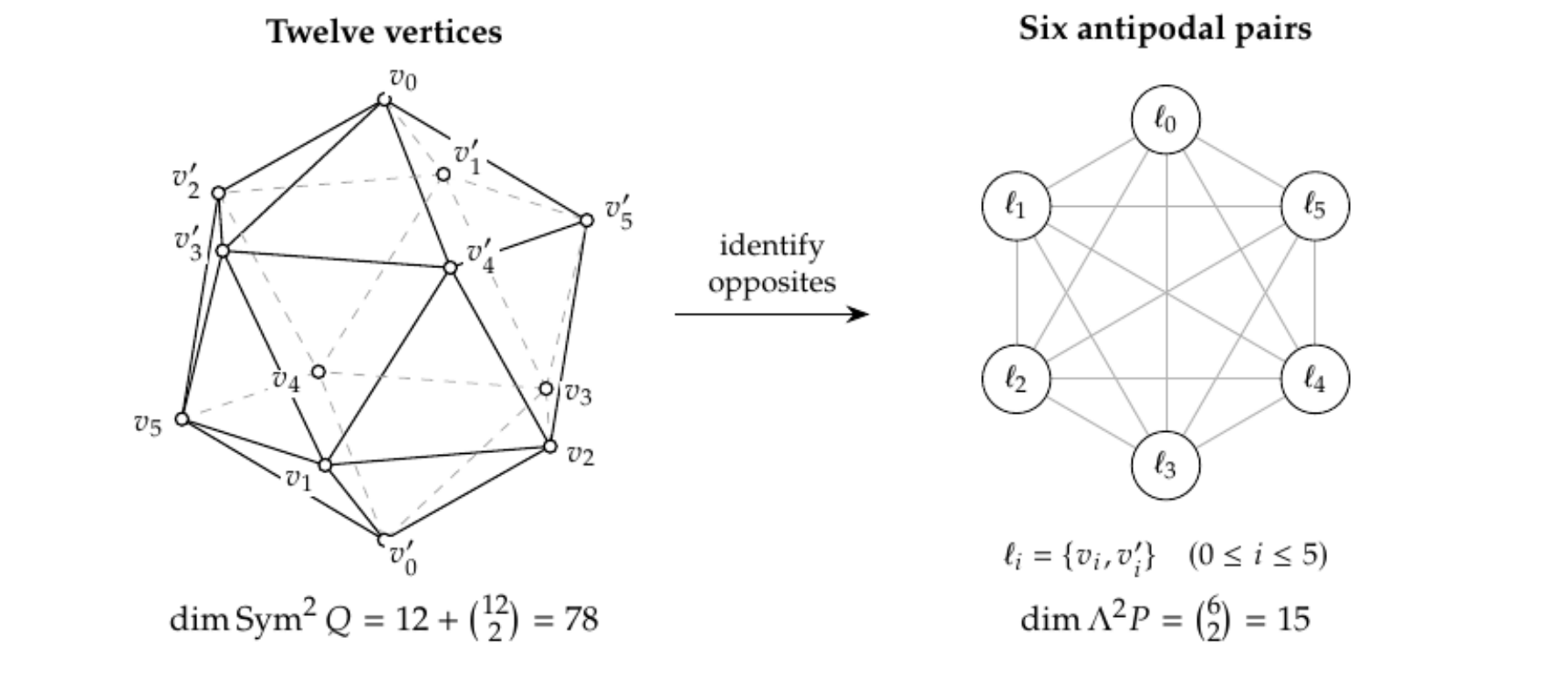}
\caption{The same twelve vertex labels will be used in Figure~\ref{fig:divisor}. Opposite vertices $v_i,v'_i$ determine the six pairs $\ell_i$; the fifteen edges of the graph on the right are the unordered pairs of axis labels. Symmetric vertex pairs permit repetition, whereas axis wedges change sign when their order is reversed.}\label{fig:vertices}
\end{figure}

Take the symmetric square of the vertex space and the exterior square of the axis space. The symmetric pairs of vertices allow a repeated vertex, whereas a wedge \(e_i\wedge e_j\) uses two distinct axis basis vectors and changes sign when their order is reversed. Hence
\[
U=\operatorname{Sym}^2Q,\quad \dim U=12+\binom{12}{2}=78,
\qquad
V=\Lambda^2P,\quad \dim V=\binom62=15.
\]
Figure~\ref{fig:vertices} labels the opposite vertex pairs; the fifteen edges of the graph on the right are the unordered pairs of axis labels.

The two square formulas applied to the character table give
\[
\chi_U=(78,6,0,3,3),\qquad
\chi_V=(15,-1,0,0,0),\qquad
\chi_{\Lambda^2U}=(3003,-21,0,3,3).
\]
For example, at a half-turn the traces on \(Q\) and \(P\) are zero and two, so the traces on $U$ and $V$ are \((0^2+12)/2=6\) and \((2^2-6)/2=-1\). It remains to compute the character of \(\Lambda^5V\).

At an involution the eigenvalues on \(V\) are \(+1\) and \(-1\). If their multiplicities are \(a,b\), then \(a+b=15\) and \(a-b=-1\), so \(a=7,b=8\). The generating polynomial becomes
\[
(1+t)^7(1-t)^8=(1-t)(1-t^2)^7,
\qquad[t^5](1-t)(1-t^2)^7=-\binom72=-21.
\]
To obtain the coefficient of \(t^5\), take \(-t\) from the first factor and the coefficient of \(t^4\) from the second, selecting \(-t^2\) from two of its seven factors. This agrees with \((6^2-78)/2=-21\) on \(\Lambda^2U\).

For an element of order three, write \(\omega=e^{2\pi i/3}\) and let \(m_0,m_1,m_2\) count its eigenvalues \(1,\omega,\omega^2\) on \(V\). The traces of the identity, \(g\), and \(g^2\) give
\[
\begin{aligned}
m_0+m_1+m_2&=15,\\
m_0+\omega m_1+\omega^2m_2&=0,\\
m_0+\omega^2m_1+\omega m_2&=0.
\end{aligned}
\]
Adding the equations and using \(1+\omega+\omega^2=0\) gives \(3m_0=15\); subtracting the last two gives \(m_1=m_2\), hence all three multiplicities are five. Each triple contributes \((1+t)(1+\omega t)(1+\omega^2t)=1+t^3\), so the full polynomial is \((1+t^3)^5\), which has no term of degree five.

At order five, Fourier inversion gives the eigenvalue multiplicities in a single formula. With \(\zeta=e^{2\pi i/5}\), the multiplicity of \(\zeta^a\) is
\[
m_a=\frac15\sum_{r=0}^{4}\chi_V(g^r)\zeta^{-ar}
     =\frac{15}{5}=3\qquad(0\leq a\leq4).
\]
Indeed, substituting the trace as a sum of eigenvalues leaves the geometric sum \(\sum_{r=0}^{4}\zeta^{(b-a)r}\), which is five when \(a=b\) and zero otherwise. All four nonidentity powers have trace zero on \(V\), including those in the other five-cycle class. The polynomial is therefore \((1+t^5)^3\), with coefficient of \(t^5\) three. At the identity, both characters equal $\binom{15}{5}=3003$. They therefore agree on all five conjugacy classes, proving
\begin{equation}\label{eq:icosahedral}
\boxed{\Lambda^2\operatorname{Sym}^2Q\cong\Lambda^5\Lambda^2P.}
\end{equation}
The same twelve vertices will later give an equivariant isomorphism for this identity through Klein's invariant polynomial, division, differentiation, and residues.

A second symmetric-square construction gives \(11628\), using the same vertex and axis modules. Adjoin \(B\) to \(Q\), obtaining \(W_{17}=Q\oplus B=\mathbf1\oplus S\), of dimension seventeen. \(W_{17}\) is also isomorphic to the zero-sum subspace of the permutation representation on the disjoint union of the twelve vertices and six axes: functions constant on each orbit have values \(a,b\), subject to \(12a+6b=0\), leaving one fixed line together with \(2B\oplus C\). Its symmetric square has dimension \(17\cdot18/2=153\). Put \(V_{19}=\mathbf1\oplus3C\), of dimension nineteen. The characters of $Q,B,C$ and the square formulas give
\[
\chi_{\operatorname{Sym}^2W_{17}}=(153,9,0,3,3),\qquad
\chi_{V_{19}}=(19,-5,1,4,4).
\]
The character of \(\Lambda^2\operatorname{Sym}^2W_{17}\) is \((11628,-36,0,3,3)\). For \(V_{19}\), the exterior polynomials at orders two, three, and five are
\[
(1+t)^7(1-t)^{12},\qquad
(1+t)(1+t^3)^6,\qquad
(1+t)^4(1+t^5)^3.
\]
The first can be rewritten as \((1-t)^5(1-t^2)^7\); its coefficient of \(t^5\) is \(-1+70-105=-36\), from degree pairs \((5,0),(3,2),(1,4)\). The other two coefficients of \(t^5\) are zero and three. Together with the dimension, this proves
\[
\Lambda^2\operatorname{Sym}^2(\mathbf1\oplus S)
\cong\Lambda^5(\mathbf1\oplus3C),\qquad
\binom{153}{2}=\binom{19}{5}=11628.
\]

\section[Faces at 1540]{Faces at $1540$}\label{sec:faces}

The twenty faces give a further permutation space, $R_{20}$. Combine
the six axes, the ten simplex edges, and two distinct copies of the face
set. For the other space, take two copies of the five letters and the
twelve vertices. Their permutation spaces are
\[
U=P\oplus E\oplus2R_{20},\qquad V=2R_5\oplus Q,
\]
with dimensions
\[
\dim U=6+10+20+20=56,\qquad \dim V=5+5+12=22.
\]
We obtain
\begin{equation}\label{eq:faces-1540}
\boxed{\Lambda^2(P\oplus E\oplus2R_{20})
       \cong\Lambda^3(2R_5\oplus Q).}
\end{equation}
Both sides have dimension $\binom{56}{2}=\binom{22}{3}=1540$.

We compute the characters by counting fixed labels. A half-turn fixes two axes
and two simplex edges, and no face, so its trace on $U$ is four. On each
five-letter copy it fixes one letter, while it fixes no vertex, giving
trace two on $V$. Hence its trace on the exterior square is
\[
\frac{4^2-56}{2}=-20.
\]
On the twenty-two labels of $V$, the half-turn has two fixed labels and ten
transpositions. Its exterior generating polynomial is
$(1+t)^2(1-t^2)^{10}$, whose coefficient of $t^3$ is also $-20$.

The characters of $U$ and $V$ are
\[
\chi_U=(56,4,5,1,1),\qquad \chi_V=(22,2,4,2,2).
\]
An element of order three fixes four labels of $V$ and permutes the
remaining eighteen in six cycles, giving
\[
[t^3](1+t)^4(1+t^3)^6=4+6=10
=\frac{5^2-5}{2}.
\]
An element in either class of order five fixes two labels of $V$ and has
four five-cycles; its polynomial $(1+t)^2(1+t^5)^4$ has no cubic term.
On $U$ the corresponding exterior-square trace is $(1^2-1)/2=0$.
The common exterior character is
\[
(1540,-20,10,0,0).
\]

Among permutation representations of dimensions $56$ and $22$, this is
the unique ordered pair of representation isomorphism classes satisfying
$\Lambda^2U\cong\Lambda^3V$, apart from the pair with both actions
trivial. Companion Section A.5 proves this by considering the traces of elements
of orders two, three, and five. Different disjoint unions of orbits can have the same permutation character;
the uniqueness is at the level of representations.

\section[Vertex triples at 24310]{Vertex triples at $24310$}\label{sec:triples}

Adding a trivial line to $\Lambda^3Q$ gives dimension
\[
221=1+\binom{12}{3}.
\]
A basis of $\Lambda^3Q$ consists of wedges $e_x\wedge e_y\wedge e_z$, one for each three-element vertex set after an ordering has been chosen. Put
\[
U=\mathbf1\oplus\Lambda^3Q.
\]
For $V$, take the five-letter space together with two copies of the axis space:
\[
V=R_5\oplus2P,\qquad \dim V=5+6+6=17.
\]
The coincidence takes the form
\begin{equation}\label{eq:triples-24310}
\boxed{\Lambda^2(\mathbf1\oplus\Lambda^3Q)
        \cong\Lambda^8(R_5\oplus2P).}
\end{equation}
Both sides have dimension $\binom{221}{2}=\binom{17}{8}=24310$.

The vertex triples also explain why $U$ admits a permutation basis.
They form five rotation orbits. Two have twenty elements: the ordinary
faces and the larger equilateral triangles whose sides join nonadjacent,
nonopposite vertices. The other three orbits have sixty elements each.
An ordered wedge can change sign under a rotation, but a rotation
preserving its underlying triple always induces an even permutation of the three vertices. Indeed, a half-turn fixes no vertex and cannot preserve a set
of odd size; the only possible nontrivial stabilizer is generated by a
three-cycle. We may therefore choose an orientation at one triple in
each orbit and transport it consistently throughout that orbit. This
gives
\[
\Lambda^3Q\cong2R_{20}\oplus3R_{60},
\]
where $R_{60}$ is the regular representation: it has one basis vector for
each rotation, with the group acting by left multiplication.
Companion Section A.6 proves the orbit count and gives the orientation argument.

The characters of $U$ and $V$ are
\[
\chi_U=(221,1,5,1,1),\qquad \chi_V=(17,5,2,2,2).
\]
For instance, a half-turn on the seventeen labels of $V$ has five fixed
labels and six transpositions. Its trace on $\Lambda^8V$ is
\[
[t^8](1+t)^5(1-t^2)^6
=5\binom62-10\binom63+\binom64=-110,
\]
which is also $(1^2-221)/2$, the exterior-square trace on $U$.
At order three the polynomial is $(1+t)^2(1+t^3)^5$; its coefficient
of $t^8$ is $\binom52=10$. At order five it is
$(1+t)^2(1+t^5)^3$, whose coefficient of $t^8$ is zero.
The exterior-square formula gives the same values, so the common character is
\[
(24310,-110,10,0,0).
\]

\section[From signs to character lines: 7140]{From signs to character lines: $7140$}\label{sec:character-lines}

The chosen orientations make $\Lambda^3Q$ a permutation representation. The construction at $7140$ uses complex lines over
finite orbits. Choosing vectors in these lines will give bases permuted
with coefficients that are cube roots of unity. A basis permuted up to
nonzero scalar multiples is called \emph{monomial}; when those scalars
are all $\pm1$, it is a \emph{signed permutation basis}.
To specify such an action, choose a line over one point of an orbit. The
subgroup fixing that point, its \emph{stabilizer}, acts on the line by a
one-dimensional character. Transporting the line around the orbit gives
one line over each point and a representation on their direct sum.

Start with the five-letter orbit. The stabilizer $A_4$ of one letter acts on the three ways to divide the
other four letters into two pairs. This action has cyclic image of order three. Choose either nontrivial
character of that image, with values $1,\omega,\omega^2$, where
$\omega=e^{2\pi i/3}$. The five lines then give a representation
isomorphic to $B$, the zero-sum part of the six-axis space.

For the second example, project the twenty face centres radially onto the unit sphere. At each image $p$, define multiplication by $i$ on $T_pS^2$ by $v\mapsto p\times v$. Each tangent plane is then a complex line.
A stabilizing rotation through $2\pi/3$ acts on it by $\omega$.
The sum of these twenty lines is
\[
T\cong A\otimes B,
\]
where $A$ is the four-dimensional zero-sum part of the five-letter space.
Companion Section A.7 verifies both line
constructions by their characters.

Put
\[
W=\mathbf1\oplus3B,\qquad
U=\Lambda^2W,\qquad V=W\oplus T.
\]
Their dimensions are $16=1+3\cdot5$, $120=\binom{16}{2}$, and
$36=16+20$. The identity becomes
\begin{equation}\label{eq:lines-7140}
\boxed{\Lambda^2(\Lambda^2W)\cong\Lambda^3(W\oplus T).}
\end{equation}
Both sides have dimension $\binom{120}{2}=\binom{36}{3}=7140$.
Their common character is
\[
(7140,-60,3,0,0).
\]
For $U=\Lambda^2(\mathbf1\oplus3B)$, companion Section A.7 proves that $W=\mathbf1\oplus3B$ is, up to isomorphism, the unique sixteen-dimensional complex representation with exterior square $U$.

The line models give bases for $W$ and $V$ in which every rotation
permutes the basis vectors and multiplies them by cube roots of unity.
For $U$, the further decomposition
\[
\Lambda^2W\cong3R_{20}\oplus3T
\]
gives such a basis as well. At an element of order three, $V$ has trace $-3$.
In a signed permutation basis, a fixed vector must have sign $+1$
under an element of order three, while a three-cycle contributes zero
to the trace. Such a trace is nonnegative. Thus $V$ admits no signed permutation basis. The bases constructed above for $U$ and $V$ use cube roots of unity.
These coefficients come from characters of the stabilizers $A_4$ and
$C_3$; $A_5$ itself has only the trivial one-dimensional character.

\section{Seven identities and their counts}\label{sec:counts}

Together with the spinorial pair \((S,E)\) at \(120\), these constructions account for all seven values. Write these representations as direct sums of $\mathbf1,A,B,C$. The character table gives \(E\cong\mathbf1\oplus A\oplus B\), and we write
\[
(a,b,c,d)=a\mathbf1\oplus bA\oplus cB\oplus dC
\]
with dimension \(a+4b+5c+6d\). Table~\ref{tab:seven} lists $U$ and $V$ in increasing order of the common dimension.

\begin{center}
\begin{minipage}{\linewidth}
\captionof{table}{Seven exterior-power identities for \(A_5\). In every row, \(\Lambda^2U\cong\Lambda^kV\), with the displayed common dimension.}\label{tab:seven}
\centering
\small
\renewcommand{\arraystretch}{1.16}
\setlength{\tabcolsep}{8pt}
\begin{tabular}{rrrrll}
\toprule
Value & \(\dim U\) & \(\dim V\) & \(k\) & \(U\) & \(V\)\\
\midrule
120 & 16 & 10 & 3 & \((0,0,2,1)\) & \((1,1,1,0)\)\\
210 & 21 & 10 & 4 & \((3,2,2,0)\) & \((1,1,1,0)\)\\
1540 & 56 & 22 & 3 & \((4,5,4,2)\) & \((3,2,1,1)\)\\
3003 & 78 & 15 & 5 & \((4,4,8,3)\) & \((0,1,1,1)\)\\
7140 & 120 & 36 & 3 & \((3,9,9,6)\) & \((1,1,5,1)\)\\
11628 & 153 & 19 & 5 & \((6,9,15,6)\) & \((1,0,0,3)\)\\
24310 & 221 & 17 & 8 & \((6,16,17,11)\) & \((3,1,2,0)\)\\
\bottomrule
\end{tabular}
\end{minipage}
\end{center}

Each representation in the table has a nontrivial action; some contain trivial summands. The construction of the spin lift will identify the first row with the restriction of Schwahn's spinorial identity.

To determine all possible pairs, we fix the two representation dimensions and the exterior degrees in each row and solve the character equations. The five complex irreducibles of $A_5$ are $\mathbf1,D_+,D_-,A,B$, with dimensions $1,3,3,4,5$. The two three-dimensional characters are exchanged by the Galois conjugation of $\mathbb Q(\sqrt5)$, and our six-dimensional module is $C=D_+\oplus D_-$. Thus $C$ is reducible over $\mathbb C$, although it has an irreducible model over $\mathbb Q$.

We count ordered pairs of representation isomorphism classes satisfying $\Lambda^2U\cong\Lambda^kV$, excluding the pair of completely trivial actions. Allowing arbitrary nonnegative multiplicities of all five irreducibles gives, in the order of Table~\ref{tab:seven},
\[
10,\quad2,\quad11,\quad19,\quad3,\quad26,\quad31.
\]
The complexifications of representations defined over $\mathbb Q$ are precisely the direct sums $a\mathbf1\oplus bA\oplus cB\oplus dC$ with $a,b,c,d\in\mathbb Z_{\geq0}$; equivalently, the multiplicities of $D_+$ and $D_-$ agree. Requiring both $U$ and $V$ to belong to this class gives
\[
4,\quad2,\quad5,\quad3,\quad1,\quad6,\quad11.
\]
At $210$, $U=\mathbf1\oplus2E$ is forced even among all complex representations, and $V$ is either $E$ or $A\oplus C$. Among representations defined over $\mathbb Q$, the fifteen-dimensional representation at $3003$ is forced to be $A\oplus B\oplus C$, and the pair at $7140$ is unique. Requiring permutation bases leaves one pair at $1540$ and two at $24310$. The larger table in companion Section A.8 also counts pairs admitting signed bases or bases permuted with arbitrary roots of unity. At $7140$, all three complex solutions have trace $-3$ on the thirty-six-dimensional representation $V$ at an element of order three. Thus none admits signed bases for both representations, and exactly one admits complex monomial bases, the construction from stabilizer characters above.

Companion Section A.4 gives the enumeration and its proof \citep{Companion2026}.

An order-two linear character gives a simple family of examples. If a finite group $G$ has a nontrivial linear character $\varepsilon$ of order two, then $\Lambda^2(n\varepsilon)$ is trivial, so a binomial identity yields an exterior-power identity by taking the representation on the other side to be trivial. The linear characters used at $7140$ belong to subgroups of $A_5$; $A_5$ itself has only the trivial linear character. Pulling our seven representations back along a surjection $G\to A_5$ gives further examples.

The triple occurrence at $3003$ allows a further test of the same equations:
\[
\binom{78}{2}=\binom{15}{5}=\binom{14}{6}=3003.
\]
For brevity, we call the exterior-power problem with representation dimensions $78$ and $15$ the $78/15$ comparison, and use analogous notation for the other two pairs. The icosahedral construction realizes the $78/15$ comparison. If complex $A_5$-modules of dimensions fifteen and fourteen satisfy $\Lambda^5U\cong\Lambda^6V$, both actions are trivial; the same conclusion holds in dimensions seventy-eight and fourteen for $\Lambda^2U\cong\Lambda^6V$. Thus exactly one of the three pairwise exterior-power identities admits a nontrivial $A_5$-realization, as proved in companion Section A.4. The $15/14$ comparison is the first Lind--Singmaster equality \citep{Lind1968,Singmaster1975}.

We return to the \(78/15\) identity to construct an \emph{intertwiner}, a linear map commuting with the group action; here the map will be invertible. The wedge bases in \eqref{eq:icosahedral} have different orbit-size distributions, as computed in companion Section B. An equivariant isomorphism must therefore use linear combinations of their vectors. Klein's twelve-point construction gives such a map after invariant splittings and identifications have been chosen.

\section[Klein's twelve points and an intertwiner at 3003]{Klein's twelve points and an intertwiner at $3003$}\label{sec:klein}

The equation of the twelve vertices allows us to construct an intertwiner through polynomial division, differentiation and residues, using small spaces of polynomials to compare the two large representations. The first step is to express the vertex and axis modules in a form on which these operations can act.

Let $H_j$ be the complex $(2j+1)$-dimensional space of degree-$j$ spherical harmonics, obtained by restricting homogeneous harmonic polynomials on $\mathbb R^3$ to the unit sphere. Through the double cover $\mathrm{SU}(2)\to\mathrm{SO}(3)$, the same representation is realized by \emph{binary forms}, homogeneous polynomials in two variables, of degree $2j$:
\[
H_j\cong\mathbb C[X,Y]_{2j}
       \cong H^0(\mathbb P^1,\mathcal O(2j)).
\]
Here $\mathbb P^1$ is the complex projective line, and $H^0$ denotes the space of global sections. The line bundle $\mathcal O(d)$ assigns a one-dimensional fiber to each point; its global sections are the degree-$d$ homogeneous forms. Even degree makes the central element $-I$ act trivially, so the binary-form representation descends to the rotation group. Restricting rotations to the icosahedron gives
\[
B\cong H_2|_{A_5},\qquad B\oplus C\cong H_5|_{A_5},\qquad
P\cong(H_0\oplus H_2)|_{A_5},\quad Q\cong(H_0\oplus H_5)|_{A_5}.
\]
These are the rotation-group models used when we decompose the symmetric and exterior powers of $P$ and $Q$. The quadratic model for $P$ comes from evaluation on the six axes. The residue construction below gives the model for $Q$. The proofs are in companion Sections F.1, F.2 and F.6 \citep{Companion2026}. All maps and decompositions in the rest of this section are taken after restriction to $A_5$.

In suitable homogeneous coordinates, Klein's invariant form is
\begin{equation}\label{eq:klein-cut}
F(X,Y)=X^{11}Y+11X^6Y^6-XY^{11}.
\end{equation}
Its twelve simple zeros form the vertex orbit $D$. We also regard $D$ as a \emph{divisor}: the sum of these twelve points, each with multiplicity one. A homogeneous form scales when the representative $(X,Y)$ of a projective point is rescaled, so its value lies in a one-dimensional fiber. Set
\[
I_j=\bigoplus_{p\in D}\mathcal O(2j)_p.
\]
This is the twelve-dimensional space of independently prescribed values in the twelve fibers. Evaluation defines a map $H_j\to I_j$; for small $j$, its image has dimension less than twelve.

A degree-$2j$ form vanishes at all twelve vertices exactly when it is divisible by $F$. Evaluation is onto as soon as $2j\geq11$, by polynomial interpolation at twelve distinct points; for our even degrees, $j\geq6$ therefore gives
\[
0\longrightarrow H_{j-6}\xrightarrow{\,F\,}H_j
\longrightarrow I_j\longrightarrow0.
\]
The sequence is exact: multiplication by $F$ is injective, and its image is precisely the kernel of the surjective evaluation map. An invariant splitting means a choice of invariant complement to this kernel; it identifies \(H_j|_{A_5}\) with \(H_{j-6}|_{A_5}\oplus I_j\).

Differentiation compares these twelve-dimensional spaces. If a local expression $f$ for $F$ changes to $uf$ under a change of frame, with $u$ nowhere zero, then
\[
d(uf)|_D=u\,df|_D,
\]
because $f$ vanishes on $D$. This is the transformation rule for a section of $\mathcal O(12)\otimes\mathcal O(-2)|_D=\mathcal O(10)|_D$, where the degree $-2$ comes from differentials on $\mathbb P^1$. The resulting section $\delta=dF|_D$ is invariant and nonzero at each vertex, since the zeros of $F$ are simple. Multiplication by $\delta$ gives
\begin{equation}\label{eq:period-five-cut}
I_j\cong I_{j+5}.
\end{equation}
The shift is five because $2j+10=2(j+5)$.

The remaining low degrees are compared by residues. The residue of a meromorphic differential at a point is the coefficient of $dz/z$ in a local coordinate $z$ vanishing there. The divisor sequence gives
\[
0\longrightarrow H_r\longrightarrow I_r
\xrightarrow{\,\beta_r\,}H_{5-r}^{*}\longrightarrow0,
\qquad 0\leq r\leq5.
\]
To describe $\beta_r$, pair prescribed values $a\in I_r$ with a binary form $b\in H_{5-r}$: \(a\,b|_D\) is a section of \(\mathcal O(10)|_D\), and division by \(F\), together with the differential, gives the residue pairing. Companion Section F.6 proves exactness directly by interpolation. In particular, at $r=0$ the quotient of $I_0=Q$ by its constant line is $H_5^*\cong H_5$, identifying $H_5|_{A_5}$ with $B\oplus C$.

This last identification also has a familiar partial-fraction form. Write $f(z)=F(z,1)$ and let $\alpha_1,\ldots,\alpha_{11}$ be its finite roots. Then
\[
\omega(z)=\sum_{i=1}^{11}r_i\frac{dz}{z-\alpha_i}
\]
has residue $-\sum r_i$ at infinity. Its numerator over the denominator $f(z)$ has degree at most ten, hence homogenizes to a form in $H_5$, and every sum-zero assignment of the twelve residues is obtained in this way.

A vertex has stabilizer $C_5$, and its fiber character determines $I_j$ by induction from this cyclic subgroup, the line construction of Section~\ref{sec:character-lines}. Frobenius reciprocity computes the multiplicity of an irreducible representation as the dimension of the eigenspace on which \(C_5\) acts by the fiber character, yielding
\[
I_0\cong H_0\oplus H_5,\qquad
I_1\cong I_4\cong H_1\oplus H_4,\qquad
I_2\cong I_3\cong H_2\oplus H_3.
\]
Choose invariant complements to split the residue sequences and equivariant isomorphisms $I_4\cong I_1$ and $I_3\cong I_2$. The splittings identify each $I_j$, for $0\leq j\leq5$, with the corresponding pair of harmonic spaces. Multiplication by $\delta$ identifies $I_j$ with $I_{j+5}$.

We apply these operations separately to the harmonic decompositions of
$\Lambda^2\operatorname{Sym}^2(H_0\oplus H_5)$ and
$\Lambda^5\Lambda^2(H_0\oplus H_2)$. For \(j\geq6\), the chosen splitting expresses \(H_j|_{A_5}\) as \(H_{j-6}|_{A_5}\oplus I_j\); periodicity replaces the quotient by \(I_r\), $r\equiv j\pmod5$, and the split residue sequence replaces it by $H_r\oplus H_{5-r}$. For example, $H_6$ becomes $H_0\oplus H_1\oplus H_4$. Repeating the process on both sides gives
\[
46H_0\oplus103H_1\oplus96H_2\oplus103H_3\oplus96H_4\oplus53H_5.
\]
Both reductions have the displayed multiplicities. Restriction of the six formal harmonic classes \(H_0,\ldots,H_5\) to \(A_5\) has rank five, with kernel generated by \(H_1+2H_2+H_3-H_4-H_5\). The difference between two such multiplicity vectors with the same restricted character is therefore an integer multiple of \((0,1,2,1,-1,-1)\); the independent calculations of the two reductions show that this integer is zero. Companion Sections F.3 and F.7 give the calculation.

Choose invariant complements and identifications between the corresponding multiplicity spaces in the two reductions. The two reductions are then isomorphisms $\Phi_L,\Phi_R$ onto the common direct sum above, giving
\begin{equation}\label{eq:residue-cut}
I_F=\Phi_R^{-1}\Phi_L:
\Lambda^2\operatorname{Sym}^2Q\xrightarrow{\sim}\Lambda^5\Lambda^2P.
\end{equation}
This gives a noncanonical intertwiner factored through division by $F$, multiplication by $\delta$, and residue maps.
\FloatBarrier
\Needspace{8\baselineskip}
\section{Returning to the first coincidence}\label{sec:spin}

\begin{figure}[!htbp]
\centering
\includegraphics[width=0.92\linewidth]{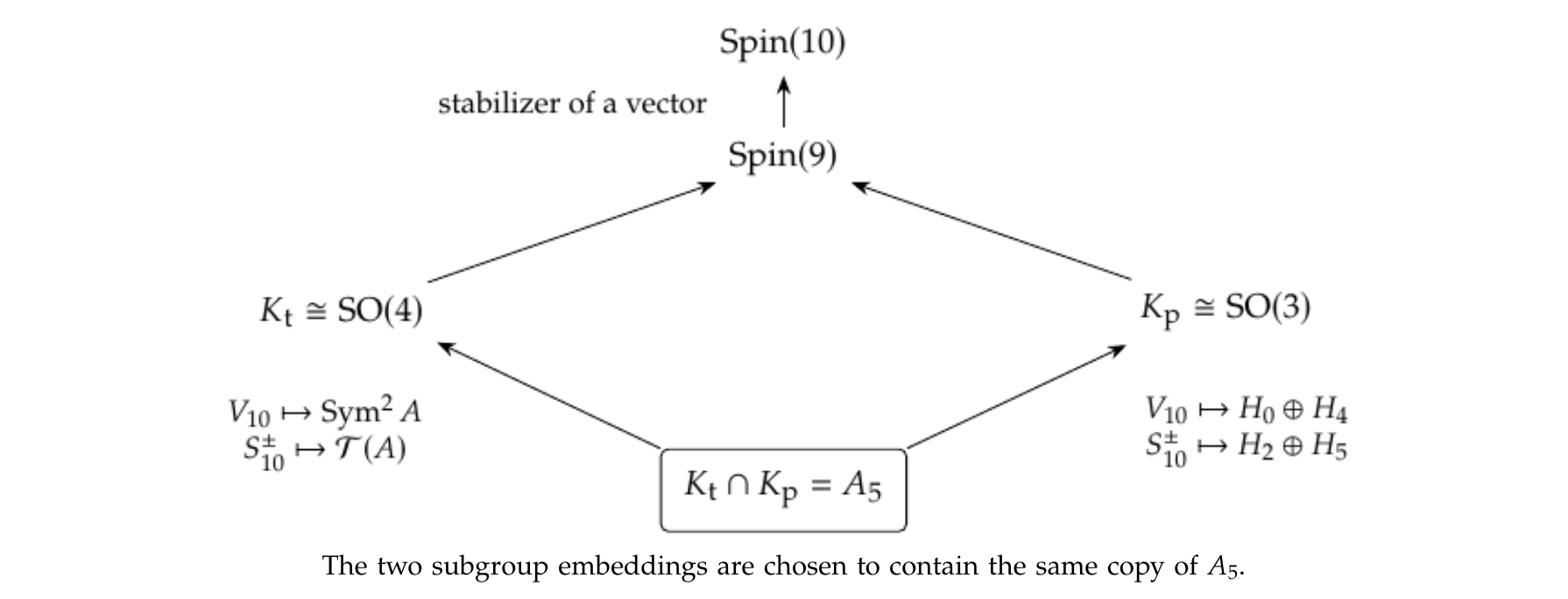}
\caption{The compatible subgroup images $K_t\cong\mathrm{SO}(4)$ and $K_p\cong\mathrm{SO}(3)$ of companion Section D.4. The first acts through the quadratic model; the second through the harmonics $H_j$. Here $\mathcal T(A)$ is the joint kernel on $A\otimes\Lambda^2A$ of wedge multiplication and metric contraction: $a\otimes(b\wedge c)$ maps to $(a\wedge b\wedge c,\langle a,b\rangle c-\langle a,c\rangle b)$. The intersection statement holds for these embeddings.}\label{fig:spin}
\end{figure}

The identity at $120$, credited to Schwahn in Baez's account, relates the ten-dimensional vector representation $V_{10}$ and either sixteen-dimensional half-spin representation $S_{10}^{\pm}$ of $\operatorname{Spin}(10)$:
\begin{equation}\label{eq:spin-classical-cut}
\Lambda^2S_{10}^{\pm}\cong\Lambda^3V_{10},
\end{equation}
and for the finite subgroup used here,
\[
V_{10}|_{A_5}\cong E,\qquad S_{10}^{\pm}|_{A_5}\cong S=2B\oplus C.
\]
The ten-dimensional restriction has an elementary quadratic model. On
\[
A_{\mathbb R}=\{(x_1,\ldots,x_5)\in\mathbb R^5:x_1+\cdots+x_5=0\},
\]
the ten functions $x_ix_j$, indexed by the simplex edges, span all quadratic forms because
\[
x_i^2=-\sum_{j\ne i}x_ix_j.
\]
The map $e_{\{i,j\}}\mapsto x_ix_j$ is equivariant. Both the edge space and the space of quadratic forms on $A$ have dimension ten, so the map is an isomorphism. The invariant metric identifies the space with $\operatorname{Sym}^2A$, giving $E\cong\operatorname{Sym}^2A$. The orthogonal symmetric-square action lifts to $\operatorname{Spin}(10)$, and the half-spin representation restricts to $S$; companion Sections D.1--D.3 give the lift and character calculations. Figure~\ref{fig:spin} shows the tensor and harmonic descriptions for the embeddings constructed in companion Section D.4.

The other six dimension triples admit no analogous nontrivial semisimple realization. The following proposition is proved as Proposition C.1 in the companion.

\begin{proposition}[Semisimple obstruction]\label{prop:semisimple}
Let $\mathfrak g$ be a finite-dimensional complex semisimple Lie algebra and let $U,V$ be finite-dimensional complex $\mathfrak g$-modules, possibly reducible. If $(\dim U,\dim V;k)$ is one of the six triples in Table~\ref{tab:seven} other than $(16,10;3)$, or is $(78,14;6)$, then $\Lambda^2U\cong\Lambda^kV$ forces the actions on both $U$ and $V$ to be trivial.
\end{proposition}

Taking $\mathfrak g=\mathfrak g_2$, $U=\mathfrak e_6|_{\mathfrak g_2}$ and $V=\mathfrak g_2$ answers Baez's exceptional-Lie-algebra question: $\Lambda^2U$ and $\Lambda^6V$ are not isomorphic, since the adjoint action on $V$ is nontrivial \citep{Baez2026Coincidences}.

The proof restricts to root copies of $\mathfrak{sl}_2$, whose modules of fixed dimension are determined by partitions of that dimension. Thus the dimension of $V$ leaves a finite list; second, fourth and sixth moments of exterior-power weights eliminate every nontrivial candidate, and the root subalgebras generate $\mathfrak g$. The $120$ identity has its spinorial source, while the remaining six have the common finite symmetry exhibited by Table~\ref{tab:seven}.

\FloatBarrier
\section{The same numbers in arithmetic}\label{sec:arithmetic}

Three of the seven values occur in a classical problem about figurate numbers, while the neighboring repetition at $3003$ belongs to the Pell--Fibonacci family.

The values $120,1540,7140$ occur as common values in the equation
\[
\binom m2=\binom n3.
\]
With $x=3(n-1)$ and $y=9m-5$, this becomes the elliptic curve
\begin{equation}\label{eq:elliptic-cut}
y^2+y=x^3-9x+20.
\end{equation}
Take the integral points
\[
P_0=(0,4),\qquad Q_0=(3,4).
\]
Starting at $P_0$, the tangent has equation $y=4-x$ and meets the curve again at $(1,3)$, so $2P_0=(1,-4)$. The line through $P_0$ and $2P_0$ is $y=4-8x$; substituting into the curve equation gives
\[
0=-x(x-1)(x-63),
\]
and hence $3P_0=(63,499)$. Returning to the binomial variables gives $n=22$, $m=56$, and therefore $\binom{56}{2}=\binom{22}{3}=1540$. The other two entries arise from $2Q_0-P_0=(27,139)$ and $-(P_0+3Q_0)=(105,1075)$, giving $120$ and $7140$. Avanesov proved that the positive integers that are both triangular and tetrahedral are exactly
\[
1,\quad10,\quad120,\quad1540,\quad7140
\]
\citep{Avanesov1967}. Thus three of the seven values occur in Avanesov's complete list.

The neighboring equality at $3003$ belongs to the Lind--Singmaster family. Cancelling factorials in
\[
\binom nk=\binom{n-1}{k+1}
\]
gives
\[
n(k+1)=(n-k)(n-k-1),
\]
which is quadratic in $n$. An integral solution $n$ therefore requires the discriminant to be a perfect square:
\[
d^2=5k^2+8k+4.
\]
Completing the square yields
\begin{equation}\label{eq:pell-cut}
(5k+4)^2-5d^2=-4.
\end{equation}
The Fibonacci and Lucas sequences satisfy the same recurrence, with initial values $(F_0,F_1)=(0,1)$ and $(L_0,L_1)=(2,1)$; their identity $L_j^2-5F_j^2=4(-1)^j$ gives solutions of this Pell equation for odd $j$. Lind proves that the positive solutions arise this way \citep[Theorem 7 and Section 5]{Lind1968}. The congruence $L_j\equiv4\pmod5$ selects $j\equiv3\pmod4$, and recovering the positive root $n=(3k+2+d)/2$ gives the normalized family, for $i\geq1$,
\[
n_i=F_{2i+2}F_{2i+3},\qquad k_i=F_{2i}F_{2i+3}.
\]
At $i=1$ this gives $(n,k)=(15,5)$, hence $\binom{15}{5}=\binom{14}{6}=3003$.

Blokhuis, Brouwer and de Weger also give families for which the difference $\binom m2-\binom n5$ is itself a triangular number, vanishing at the parameters yielding $3003$ and $11628$ \citep[Section 3.2]{BlokhuisBrouwerdeWeger2017}; companion Section E.2 gives their factorizations.

\FloatBarrier
\section{The twelve points once more}\label{sec:cusps}

Klein's twelve-point divisor also occurs in the theory of modular curves. Let $\Gamma(5)$ be the subgroup of matrices in $\mathrm{SL}_2(\mathbb Z)$ congruent to the identity modulo five. It acts on the upper half-plane $\mathfrak H=\{\tau:\operatorname{Im}\tau>0\}$ by $\gamma\tau=(a\tau+b)/(c\tau+d)$. The quotient $\Gamma(5)\backslash\mathfrak H$ is not compact. Its missing ends are filled by adding \emph{cusps}: one point for each $\Gamma(5)$-orbit in $\mathbb Q\cup\{\infty\}$ on the boundary of $\mathfrak H$. Here there are twelve. The resulting compact curve $X(5)$ is a projective line, and its cusps form the same $A_5$-set as the twelve vertices of the icosahedron \citep[Sections 3 and 6]{Duke2005}.

At the cusp represented by infinity, a local coordinate is $u=e^{2\pi i\tau/5}$: as $\operatorname{Im}\tau\to\infty$, we have $u\to0$, and adjoining $u=0$ fills the end. At another cusp, first move its representative to infinity by a modular transformation and use the same coordinate there.

A \emph{weight-two modular form} is a holomorphic function $f$ on $\mathfrak H$ satisfying
\[
f(\gamma\tau)=(c\tau+d)^2f(\tau)\qquad(\gamma\in\Gamma(5)),
\]
and holomorphic at every cusp. The latter condition means that, after moving the cusp to infinity and transforming $f$ by the weight-two rule, its expansion is $a_0+a_1u+\cdots$ with no negative powers. Since
\[
d(\gamma\tau)=\frac{d\tau}{(c\tau+d)^2},
\]
the differential $f(\tau)d\tau$ is invariant and descends to the quotient. Near a cusp,
\[
d\tau=\frac{5}{2\pi i}\frac{du}{u},
\]
so the descended differential has at most a simple pole, with residue $5a_0/(2\pi i)$. Conversely, a differential with at most simple poles at the cusps pulls back to $f(\tau)d\tau$ with this transformation law and these holomorphy conditions.

Let $M_2(\Gamma(5))$ denote the space of these modular forms. Residues give the equivariant map
\[
M_2(\Gamma(5))\longrightarrow
\left\{(r_p)_{p\in D}:\sum_{p\in D}r_p=0\right\}.
\]
If all residues vanish, the differential has no poles and is zero on $X(5)\cong\mathbb P^1$. Every sum-zero residue vector, on the other hand, is realized by the partial fractions in Section~\ref{sec:klein}. This proves both bijectivity and the dimension $12-1=11$, giving
\[
M_2(\Gamma(5))\cong B\oplus C,\qquad
Q\cong\mathbf1\oplus M_2(\Gamma(5)),\qquad
S\cong M_2(\Gamma(5))\oplus B.
\]
Under a suitable projective coordinate the cusps are the zeros of \(F\), so the eleven-dimensional zero-sum vertex representation is simultaneously the residue representation of the modular curve. Here $M_2$ denotes all holomorphic modular forms of weight two. A \emph{cusp form} has constant term $a_0=0$ at every cusp, so its differential has no poles. There are no nonzero holomorphic differentials on $\mathbb P^1$; hence the space of cusp forms is zero.

\begin{figure}[!htbp]
\centering
\includegraphics[width=0.92\linewidth]{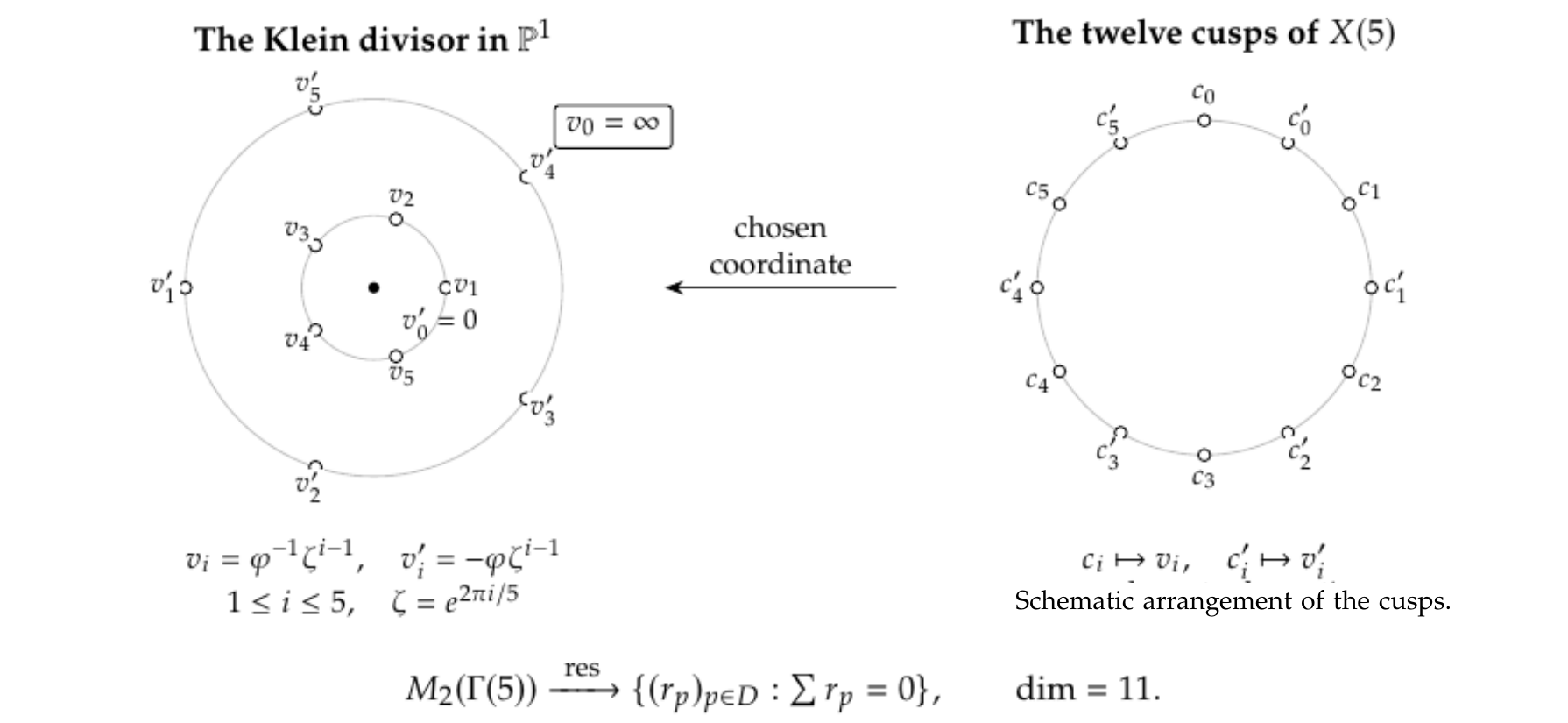}
\caption{Klein's twelve-point divisor and the twelve cusps of $X(5)$ under a chosen projective coordinate. The common $A_5$-set gives both the vertex module $Q$ and the residue description of $M_2(\Gamma(5))$.}\label{fig:divisor}
\end{figure}
\FloatBarrier

\section{A Conway epilogue}\label{sec:conway}

Companion Section H.1 exhibits a subgroup \(A_5\subset M_{24}\) with two orbits of twelve letters. Each stabilizer has order five, so each orbit is an $A_5/C_5$ coset space and carries the vertex representation $Q$. Writing $\mathfrak a$ for the complexification of the twenty-four-dimensional representation used in the Conway construction, we obtain
\[
\mathfrak a|_{A_5}\cong Q\oplus Q.
\]
The conformal construction of Duncan--Mack-Crane has weight-one space $\Lambda^2\mathfrak a$ \citep{DuncanMackCrane2015}. With the trivial action on $\mathbb C^2$, the decomposition of the exterior square of a tensor product gives
\[
\Lambda^2(Q\otimes\mathbb C^2)
\cong3\Lambda^2Q\oplus\operatorname{Sym}^2Q,
\]
where the multiplicities are $\dim\operatorname{Sym}^2\mathbb C^2=3$ and $\dim\Lambda^2\mathbb C^2=1$. The $78$-dimensional representation used in the $78/15$ equality at $3003$ consequently occurs in this weight-one space. The subgroup calculation and the tensor decomposition are given in companion Section H \citep{Companion2026}.

\section{The seven entries again}\label{sec:ending}

The seven identities in Table~\ref{tab:seven} answer Baez's question using representations of a single finite group. The first restricts the spinorial identity credited to Schwahn; the others use simplex edges, vertices, axes, faces, and stabilizer characters. The enumeration singles out the permutation pair at $1540$ and the complex monomial pair at $7140$, where signs alone cannot suffice.

At $3003$, only the $78/15$ comparison admits a nontrivial $A_5$-action. Klein's divisor gives a factored intertwiner: division by $F$ lowers degree, multiplication by $dF|_D$ gives $I_j\cong I_{j+5}$, and residues compare the low degrees. The construction depends on chosen invariant splittings and multiplicity-space identifications. Its twelve points reappear as the cusps of $X(5)$ and, in two copies, in the Conway construction.

Taking dimensions in the seven representation identities recovers
\[
120,\quad210,\quad1540,\quad3003,\quad7140,\quad11628,\quad24310,
\]
the seven values in Singmaster's short addition, all accompanied by the same sixty symmetries.

\FloatBarrier
\clearpage
\section*{Acknowledgments}

We are grateful to John Baez for posing the question that motivated this article and for reporting Paul Schwahn's spinorial explanation of the first coincidence. E.L. acknowledges support from CONAHCYT grant CB-2017-2018-A1-S-30345, the Simons Foundation grant SFI-MPS-T-Institutes-00007697, and the Ministry of Education and Science of the Republic of Bulgaria under contract DO1-239/10.12.2024. He thanks CINVESTAV for his sabbatical leave, and the International Center for Mathematical Sciences--Sofia at the Institute of Mathematics and Informatics, Bulgarian Academy of Sciences, and the Institute for the Mathematical Sciences of the Americas at the University of Miami for their hospitality and support.

\medskip
ChatGPT and Claude were used to assist with verifying mathematical proofs, checking calculations, and improving the English exposition. The authors take responsibility for the results and their presentation.

\end{document}